# Application Potential of a Hybrid Ground Source Heat Pump Array for the UC Berkeley Campus Business and Law Node Energy System: A Preliminary Study

Kecheng Chen[1], Kenichi Soga[1], Patrick Dobson[2] and Peter Nico[2]

1. University of California, Berkeley, CA, USA
2. Lawrence Berkeley National Laboratory, Berkeley, CA, USA

kecheng_chen@berkeley.edu

## ABSTRACT

The current plan divides the UC Berkeley (UCB) campus energy system into five nodes, where the Business and Law node was studied because of an open field site for borehole installation. The Pacific Northwest National Laboratory's Commercial Prototype Building Models were used to estimate heating and cooling load requirements for UCB campus building types by considering model characteristics (for example, high base load from hospitals, high DHW in hotels) corresponding to the ASHRAE Standard 90.1-2013. Unscaled load profiles were created from the EnergyPlus building energy simulation and scaled with monitored peak load and annual energy use to generate the target node's hourly heating and cooling load profiles. An optimization problem was solved to design a hybrid GSHP system, where the objective function is the lifetime total cost of the system, and the optimization variables are the portion of heating and cooling loads covered by the GSHP system. Modelica models for air source and ground source heat pump systems were built for detailed case studies based on optimization results. In the Modelica model, the demand side is connected to the radiators in the building to transfer heat, and the source side is connected to GSHP, ASHP, or other heating and cooling facilities. The results demonstrate that an appropriate hybrid GSHP system can help reduce both borehole numbers and electricity consumption for the UCB campus site.

## 1. INTRODUCTION

Because of the increasingly deteriorating state of the existing campus energy delivery system, UCB recognized the need for a holistic and long-term study of the future of its campus energy delivery system. Under the "UCB Energy Delivery Options Analysis" project, UCB contracted an engineering consultancy company, Arup, to study the best method of delivering heat and power to the campus (Arup, 2015). Considering factors including general use and occupancy profile, proximity, and load densities, the nodal approach was determined for effective energy delivery. It was proposed to divide the campus into five nodes, as shown in Figure 1.

As a reliable, high-efficient, and carbon-free energy source, shallow geothermal or ground source heat pump (GSHP) has got attention from the campus to achieve the "Clean Energy Campus" goal (Sun et al., 2021; Chen et al., 2021; Chen et al., 2021). In unbalanced load conditions, like the UCB campus, whose cooling load is much higher than the heating load, it commonly complains that the capacity of the GSHP became deficient in maintaining a comfortable indoor air temperature after a few years of operation (Xu et al., 2021). The soil thermal imbalance mainly causes this. A hybrid GSHP system is the universal alternative to deal with thermal imbalance (Qi et al., 2014). It also helps to make GSHP more cost-effective by reducing the high first cost of borehole installation (Kavanaugh, 1998). Considering the carbon-free energy goal, the whole electrical hybrid system combining GSHP and the air source heat pump (ASHP) would be a favorable option for the UCB campus (You et al., 2021; Martins and Bourne, 2022). As shown in Figure 2, UCB plans to have an open field in the Business and Law node (the southeast corner of the campus) for borehole installation and campus geothermal potential study.

This work focuses on the hybrid GSHP study in the Business and Law node, which can extend to the whole campus in the future. The data for each building in the node provided by Arup is shown in Table 1 (Arup, 2015). Peak heating load was calculated from peak steam intensity, and peak cooling load was calculated from total chiller capacity. In this study, unscaled annual space heating/cooling and domestic hot water (DHW) load profiles were developed for each building using EnergyPlus simulations with Commercial Prototype Building Model (CPBM) (Crawley et al., 2001; Goel et al., 2014). The unscaled load profiles were scaled to match the peak load and annual energy consumption to ensure the actual demands were considered. This will be described in section 2.1. With the scaled load profiles, an optimization problem was formed to minimize the system's 20-year total cost, including installation and operation costs, and to design the hybrid GSHP/ASHP system. This will be described in section 2.2. Modelica-based dynamic GSHP/ASHP system was built to study the transient state of the system in detail. This will be described in section 2.3. With the optimal design results and the dynamic simulation model, case studies were done by changing the number of boreholes and the capacity of GSHP. The results will be described in section 3.

Chen et al.

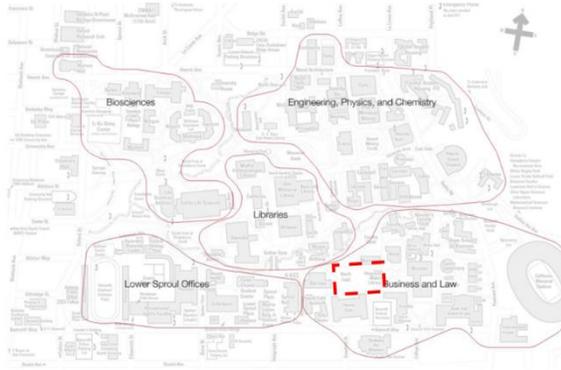

**Figure 1: Campus zoning – nodal energy delivery approach (open field for borehole installation shown in red rectangular)**

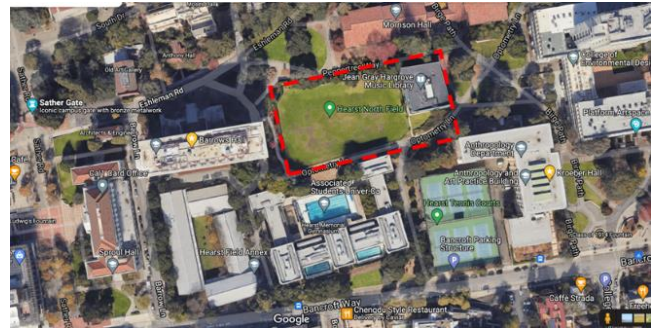

**Figure 2: Satellite image of the Business and Law node (open field for borehole installation shown in red rectangular)**

**Table 1: Building data of the Business and Law node**

| Building name | Building Area (sqft) | Annual Electricity Use (kWh/y) | Annual Steam Use (kBtu/y) | Peak Heating (kBtu/h) | Peak Cooling (tons) |
|---|---|---|---|---|---|
| Barrows Hall | 193202 | 1260397 | 4967095 | 1763 | 193 |
| Boalt Hall | 277521 | 1833318 | 7134880 | 3429 | 278 |
| California Memorial Stadium | 261091 | 2692449 | 1876213 | 5397 | 261 |
| Calvin Laboratory | 33167 | 370341 | 4199765 | 656 | 85 |
| Cheit Hall | 32120 | 604834.8 | 825783.8 | 864 | 32 |
| Faculty Club | 35483 | 337170 | 3167710 | 802 | 35 |
| Gerson Bakar Faculty Building | 106295 | - | - | 2747 | 106 |
| Hargrove Music Library | 29063 | 341660 | 827515 | 563 | 29 |
| Hearst Field Annex | 34825 | 81065 | - | - | 52 |
| Hearst Memorial Gymnasium | 124197 | 602605 | 20178687 | 2728 | 124 |
| Hertz Memorial Concert Hall | 30160 | 327029 | 1710967 | 514 | 30 |
| International House | 185200 | 1432521 | 16533547 | 4786 | 130 |
| Kroeber Hall | 117814 | 682886 | 4351656 | 780 | 118 |
| Minor Hall | 46204 | 58723 | 48477142 | 669 | 69 |
| Minor Hall Addition | 55516 | 1249244 | 4217821 | 1435 | 500 |
| Morrison Hall | 40706 | 180419 | 2309238 | 953 | 41 |
| Simon Hall | 36180 | 239007 | 1362496 | 935 | 36 |
| Simpson Student-Athlete High-Performance Center | 139682 | 2246028 | 2606521 | 677 | 140 |
| Sproul Hall | 110919 | 910249 | 3355426 | 1370 | 111 |
| Walter A. Haas Student Services Building | 95712 | 1802302 | 2460692 | 2473 | 200 |
| Women's Faculty Club | 18474 | 98954 | 1649248 | 418 | 18 |
| Wurster Hall | 222434 | 1750749 | 3615932 | 898 | 222 |

## 2. MODEL ESTABLISHMENT

### 2.1 Commercial Prototype Building Model and EnergyPlus Simulation

As part of DOE's support of ANSI/ASHRAE/IES Standard 90.1 and IECC, PNNL developed a suite of prototype buildings, which cover 75% of the commercial building floor area in the United States for new construction, including both commercial buildings and mid- to high-rise residential buildings across all US climate zones (Goel et al., 2014). For the UCB project, Arup created a mapping table between its CPBM and the UCB campus building types by considering model characteristics (for example, high base load from hospitals and high DHW in hotels) corresponding to the ASHRAE Standard 90.1-2013 (Arup, 2015). The mapping relationship is given in Table 2.

CPBM provides EnergyPlus model input files that can be used to conduct building simulations. EnergyPlus is a whole building energy simulation program that models energy consumption—for heating, cooling, ventilation, lighting, and plug and process loads—and water use in buildings (Crawley et al., 2001). Based on the mapping relationship given in Table 2 and Berkeley's typical meteorological year (TMY) weather data collected from ResStock (Wilson, 2017), unscaled annual space heating/cooling and domestic hot water (DHW) load profiles are developed for each building. Assuming an additional facility is needed to lift the water temperature for DHW use, the space heating and DHW loads can be combined as the heating load.

As shown in Table 3, the parameters proposed by Arup were used to generate the end-use customer's benchmark peak load and annual energy consumption from Table 1 (Arup, 2015). The benchmark peak heating/cooling load for the Business and Law node was calculated by adding peak heating/cooling for each building in the node and multiplying the node load diversity coefficient. The original annual energy consumption data were counted separately based on the steam use and electricity use. Building heating system efficiency was used to convert the steam used to the annual heating energy consumption. The average chiller coefficient of performance (COP) for the node was used to convert the electricity usage to the annual cooling energy consumption. As the chiller COP changes concerning the ambient wet bulb temperature, a more reasonable seasonal COP was used.

The following optimization function was formed to scale the unscaled load profiles to the benchmark peak load and annual energy consumption.

$$\min_{k>0,b} \left\| \frac{Sum(k*L_{unscaled}+b) - E_{annual}}{E_{annual}} \right\|_2^2 + \left\| \frac{Max(k*L_{unscaled}+b) - P}{P} \right\|_2^2 \qquad (1)$$

Where $L_{unscaled}$ is the unscaled heating or cooling load profile for the Business and Law Node, $E_{annual}$ is the annual heating or cooling energy consumption, and $P$ is the peak heating or cooling load. Python scipy package was used to find the optimal solution for solving this unconstrained multivariate scalar minimization function (Virtanen et al., 2020). The resulting scaled profiles are shown in Figure 3. After checking, the annual energy consumption and peak load calculated from scaled load profiles were well matched with the metered data.

**Table 2: Mapping between the UCB campus building type and the CPBM**

| UCB campus building type | CPBM |
|---|---|
| Admin | Large office |
| Athletic | Large hotel |
| Datacenter | large hotel |
| Auditorium | Large office |
| Bioscience | Hospital |

Chen et al.

| | |
|---|---|
| Chemistry | Hospital |
| Engineering | Hospital |
| Classroom | Secondary school |
| Residence | Mid-rise apartment |
| Mixed-use | Average of the hospital, secondary school, and large office |

**Table 3: Parameters for benchmark load data generation**

| Parameter | Value |
|---|---|
| Existing building-level steam loss | 15% |
| Gas efficiency | 80% |
| Absorption chiller efficiency | 0.4 |
| Existing electric chiller efficiency | 1.25kW/ton |
| New electric chiller efficiency | 0.5kW/ton |
| Building heating system efficiency | 95% |

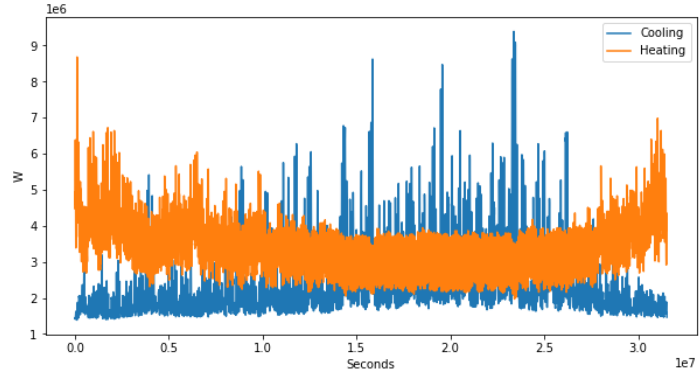

**Figure 3: Scaled heating and cooling load profiles of the Business and Law node**

**2.2 Optimization-based Hybrid Ground Source Heat Pump Design**

As shown in Eq (2), an optimization problem was solved to design the hybrid GSHP, where the objective function is the net present value of the 20-year total cost of the hybrid system (Alavy et al., 2013).

$$\min_{0<\alpha<1, 0<\beta<1} \sum_{t=0}^{20} \frac{CF_t}{(1+IR)^t} \quad (2)$$

Where $CF_t$ is the cash flow at time t and IR is the interest rate. The cost calculation parameters are shown in Table 4. To find an optimal GSHP capacity in the hybrid system, shave factors, $\alpha$ for cooling and $\beta$ for heating were introduced. $\alpha(\beta)$ is the portion of the cooling(heating) load covered by GSHP ($\alpha$ in Figure 4). The rest of the loads are covered by ASHP. The greedy algorithm was used to solve the current optimization problem by searching $\alpha$ from 0 to 1 with a step size of 0.01. Because the Business and Law Node was cooling dominated, as $\alpha$ was varied from 0 to 1, a value of $\beta$ was found, which corresponds to the portion of heating demand that would be met by a GSHP that meets a portion of the cooling demand corresponding to $\alpha$. The COP of ASHP and GSHP were assumed to be 2.5 and 3.5 in the heating mode and 4.5 and 5.5 in the cooling mode, respectively (Chen et al., 2021). Calculation of the required borefield lengths for cooling and heating was based on Eqs (3) and (4) presented by Kavanaugh and Rafferty (Kavanaugh and Rafferty, 2014).

$$L_c = \frac{\left(q_a R_{ga} + \left(C_{fc} q_{lc}\right)\left(R_b + PLF_m R_{gm} + R_{gd} F_{sc}\right)\right)}{t_g - \frac{t_{wi} + t_{wo}}{2} - t_p} \quad (3)$$

$$L_h = \frac{\left(q_a R_{ga} + \left(C_{fh} q_{lh}\right)\left(R_b + PLF_m R_{gm} + R_{gd} F_{sc}\right)\right)}{t_g - \frac{t_{wi} + t_{wo}}{2} - t_p} \quad (4)$$

$$q_a = \frac{C_{fc} q_{lc} EFLH_c + C_{fh} q_{lh} EFLH_h}{8760 \text{ hours}} \quad (5)$$

Where $q_a$ can be calculated from Eq (5), and other terms are interpreted in Table 5. Assumptions for the design are shown as the following.

- 1.25" HDPE single vertical U-pipe
- The borehole diameter is 0.127 m
- The static water table at 5 m below surface
- The period lengths of annual, monthly, and daily heat pulses are twenty-year (7300 days), one month (30 days), and six hours (0.25 days)
- For the ground, the undisturbed temperature is 18 °C, the thermal conductivity is 2.42 $W/(mK)$ and diffusivity is 0.08 $m^2/day$
- Borehole thermal conductivity is 1.4 $W/(mK)$
- For the heat pump, the entering liquid temperature is 25 °C (cooling) / 8°C (heating), and the leaving liquid temperature is 30 °C (cooling) / 3 °C (heating)

Figure 5 shows the initial cost, the operation cost, and the total cost as functions of α. As $\alpha$ increases, the GSHP capacity, and thus its drilling and installation costs increase, but the operation cost decreases, as the Business and Law Node will be relying less on the conventional systems. The optimal total cost of the GSHP/ASHP combination system is $19.1 million for a 20-year operation, where 88% of the cooling load and 98% of the heating load can be covered by GSHP with 231*200m boreholes. The optimal base load that should be covered by GSHP is shown in Figure 6.

**Table 4: Parameters for cost calculation**

| Parameter | Cost/Value |
|---|---|
| The ground heat exchanger (installation and materials) | $65.5/meter |
| Air-to-water heat pump | $80/kW of heat pump design capacity |
| Industrial electricity | $0.08/kWh |
| Interest rate (IR) | 8% |
| Inflation rate | 4% |

**Table 5: Interpretation of terms in Kavanaugh and Rafferty's design formula**

| Term | Interpretation |
|---|---|
| $F_{sc}$ | Short-circuit heat loss factor |
| $L_c$ and $L_h$ | Required ground-loop lengths to meet the shaved cooling and heating loads |
| $F_m$ | Part-load factor during the design month |
| $q_a$ | Net annual average heat transfer to the ground |
| $q_{lc}$ and $q_{lh}$ | Building design cooling and heating block loads |
| $R_{gd}, R_{gm}$ and $R_{ga}$ | Daily, monthly, and annual effective thermal resistance of the ground |

Chen et al.

| $R_b$ | The thermal resistance of the borehole |
|---|---|
| $t_g$ | Undisturbed ground temperature |
| $t_p$ | Temperature penalty |
| $t_{wi}$ and $t_{wo}$ | Water temperature at heat pump inlet and outlet |
| $C_{fc}$ and $C_{fh}$ | Correction factors depending on COP in the cooling and heating mode |
| $EFLH_c$ and $EFLH_h$ | Annual equivalent full-load cooling and heating hours |

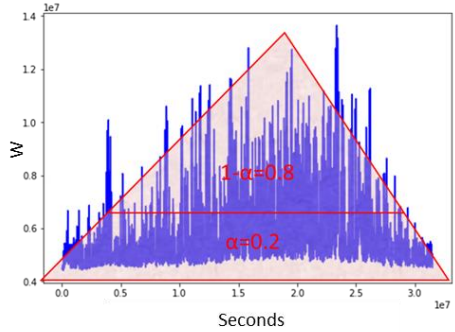

Figure 4: Explanation of $\alpha$ ($\alpha$ represents the ratio between the area of the base cooling load and the area of the total cooling load)

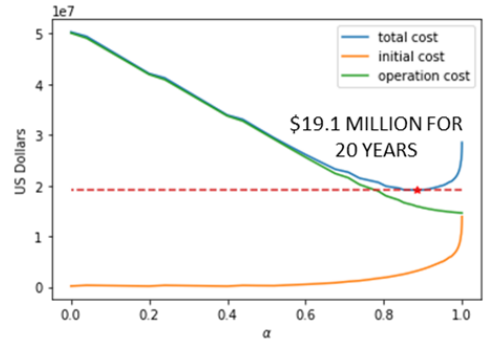

Figure 5: The initial cost, the operation cost, and the total cost as functions of $\alpha$

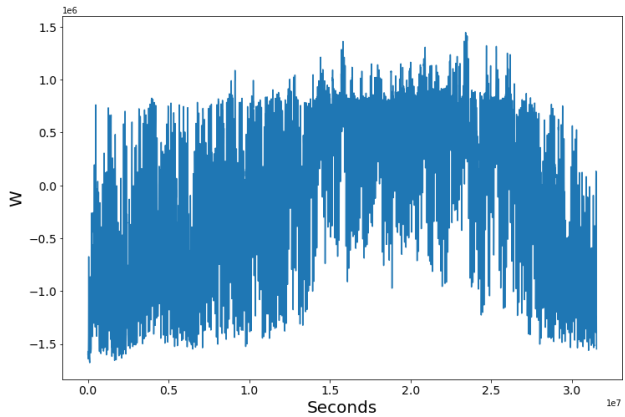

Figure 6: Optimal base load covered by GSHP

**2.3 Modelica-based Hybrid Ground Source Heat Pump Modeling**

The structure of the proposed Modelica model is shown in Figure 7. In the model, the demand side is connected to the radiators in the building to transfer heat. The source side is connected to GSHP, ASHP, or the combination. As shown in Figure 8, the model was built based on the Modelica Buildings Library (Wetter et al., 2014). It mainly consists of five parts, including building load, outdoor weather, control loop, ASHP, and GSHP. For the building load part, the net heating and cooling loads of the Business and Law node in Figure 3 were applied to a mixing volume whose size equals the total floor area of buildings multiplied by the floor height. The mixing volume was connected to an air radiator which models both convective and radiative heat transfer. The weather data in section 2.1 was used for the outdoor weather part.

Two control strategies were proposed for the control loop to keep the indoor temperature from 21 °C to 24 °C. The first control strategy is that when the indoor temperature is higher than 26 °C (cooling mode) or lower than 19 °C (heating mode), the water pumps at the demand and source sides of the heat pump would be activated. In cooling mode, the water pumps keep running until the indoor temperature reaches 21 °C. In heating mode, the water pumps keep running until the indoor temperature reaches 24 °C. When the water pumps stop, the system is in sleep mode. The nominal flow rate of the water pump was set to make the water loop have a 5 °C temperature difference when experiencing the average load. The second control strategy is that when the system is not in sleep mode, and the transient flow rate of the water pump reaches the heat pump's minimal operation flow rate, the heat pump would be activated; otherwise, it stops. The production water temperature from the heat pump to the air radiator was set to 50 °C in the heating mode and set to 10 °C in the cooling mode. ASHP part would be activated when the operation water temperature of the ground loop was higher than 30 °C or lower than 5 °C.

For both ASHP and GSHP parts, the reversible heat pump was used because it can be operated either in heating or cooling mode by reversing the flow of refrigerant from the compressor through the condenser and evaporation coils. The performance data of GSHP was taken from Trane EXW240, and that of ASHP was taken from Trane PICCO270 (ASHRAE Handbook, 2012). HP performances were incorporated into the model through the equation fit method by fitting Eqs (6) and (7) to the performance data and finding the values of the performance coefficients.

$$\dot{Q} = \left(\alpha_1 + \alpha_2\left(\frac{T_{loa}}{T_{Refloa}}\right) + \alpha_3\left(\frac{T_{sou}}{T_{Refsou}}\right) + \alpha_4\left(\frac{\dot{m}_{loa}}{\dot{m}_{Refloa}}\right) + \alpha_5\left(\frac{\dot{m}_{sou}}{\dot{m}_{Refsou}}\right)\right)\dot{Q}_{Ref} \qquad (6)$$

$$P = \left(\beta_1 + \beta_2\left(\frac{T_{loa}}{T_{Refloa}}\right) + \beta_3\left(\frac{T_{sou}}{T_{Refsou}}\right) + \beta_4\left(\frac{\dot{m}_{loa}}{\dot{m}_{Refloa}}\right) + \beta_5\left(\frac{\dot{m}_{sou}}{\dot{m}_{Refsou}}\right)\right)P_{Ref} \qquad (7)$$

Where $\dot{Q}$ is the heat transfer from the load side, $P$ is the power consumption, $T$ is the liquid temperature, $\dot{m}$ is the flow rate, the subscript $loa$ means the load side, the subscript $sou$ means the source side, the subscript $Ref$ means the reference design value, and $\alpha_n$ and $\beta_n$ (n=1,2,3,4,5) are the performance coefficients. For the ASHP part, the outdoor climate determined the entering air temperature. For the GSHP part, the borehole was modeled using an axial discretization and a resistance-capacitance network for the internal thermal resistances between the individual pipes and between each pipe and the borehole wall. The ground properties that appeared in section 2.2 was used. The ground response was modeled through the spatiotemporal superposition of the modified g-functions. The final measured electricity consumption includes the consumptions from the water pumps and the heat pumps.

A case was done to test the model, where the optimal base load covered by GSHP and the required borehole length calculated from section 2.2 were used to set the building load and the borehole configuration. GSHP-only case was studied. Figure 9 shows the annual measured indoor temperature. The indoor temperature could stay inside 21-24 °C for most of the time and matched the proposed control strategies well. Figure 10 shows the temperature of the water supplied to and returned from the ground. The temperature difference matched 5 °C assumption well. The correctness of the model was verified.

## 3. RESULTS

The scaled heating and cooling load profiles for the Business and Law Node shown in Figure 3 were used to set the building load. Case studies were done by changing the number of boreholes and the capacity of GSHP. The resulting annual electricity consumption of the hybrid systems is shown in Figure 6. GSHP-only system with 1060 boreholes can decrease 28% electricity consumption compared to the ASHP-only system. However, the installation cost of the borehole is high, and the borehole number should be further reduced. GSHP-only system with 616 boreholes can have a 0.06% electricity consumption increase and a 42% borehole number decrease compared to the GSHP-only system with 1060 boreholes. The GSHP-ASHP hybrid system was also studied with 616 boreholes, where the role of ASHP was to take a part of loads from the ground. It was found that the GSHP-ASHP hybrid system with 616 boreholes can have a 0.02% electricity consumption decrease compared to the GSHP-only system with 616 boreholes. So the hybrid configuration can help decrease electricity consumption by improving the efficiency of GSHP. GSHP-only system with 465 or fewer boreholes failed to work because the operation water temperature of the ground loop was beyond the acceptable range. Another common electrical equipment, the electric heater, was also studied in the hybrid configuration. GSHP-electric heater hybrid system with 465 boreholes can have a 19% electricity consumption increase compared to the ASHP-only system. It is because the heat pump uses electricity to transfer the heat and the electric heater uses electricity to produce the heat. The heat pump has lower electricity consumption for the same amount of heat. GSHP-ASHP hybrid system with 231 and 150 boreholes can decrease electricity consumption by 14% and 0.09% compared to the ASHP-only system.

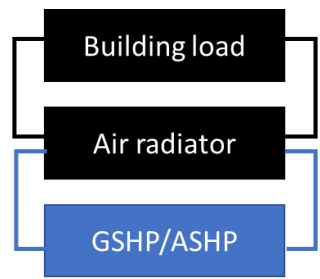

Figure 7: The structure of the proposed Modelica model

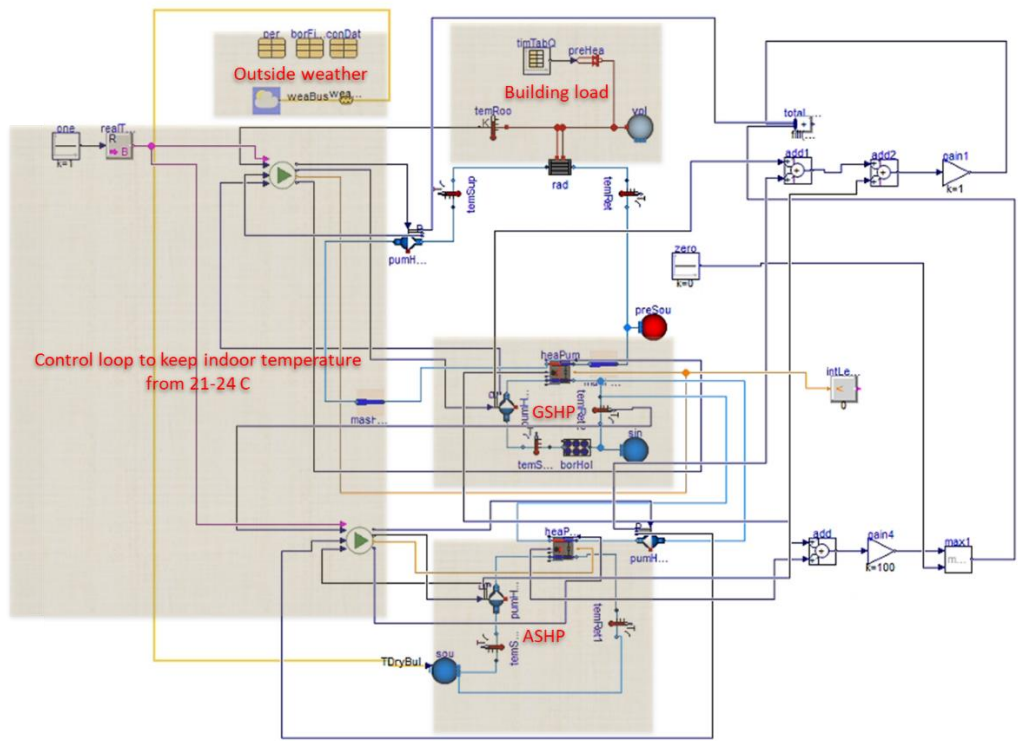

Figure 8: Dynamic model of the hybrid GSHP/ASHP system

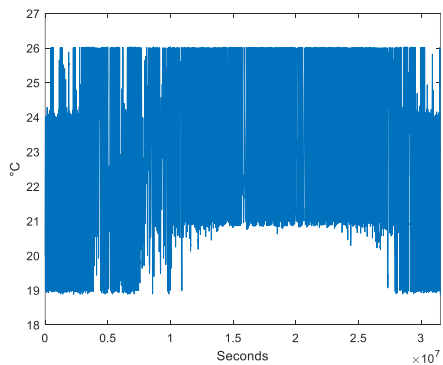

Figure 9: The annual measured indoor temperature

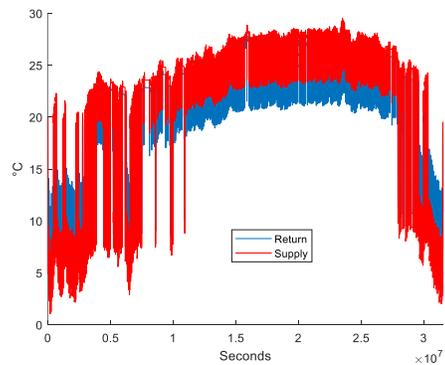

Figure 10: The annual temperature of the water supplied to and returned from the ground

Table 6: Case study results

| Case | Option | Borehole number (each 200m) | Electricity consumption (kWh) |
|---|---|---|---|
| 1 | ASHP only | 0 | 8.9637E+05 |
| 2 | GSHP only | 1060 | 6.4791E+05 |
| 3 | GSHP only | 616 | 6.4828E+05 |
| 4 | GSHP + ASHP | 616 | 6.4815E+05 |
| 5 | GSHP only | 465 | Fail to work |

| 6 | GSHP + Electric heater | 465 | 1.0636E+06 |
| 7 | GSHP only | 231 | Fail to work |
| 8 | GSHP + ASHP | 231 | 7.6848E+05 |
| 9 | GSHP only | 150 | Fail to work |
| 10 | GSHP + ASHP | 150 | 8.8876E+05 |

## 4. CONCLUSION

The Business and Law node's hourly heating and cooling load profiles were generated by scaling the prototypical building model outputs. The annual energy consumption and peak load calculated from scaled load profiles were well-matched with the metered data. Hybrid GSHP was designed for the node with an optimization-based method. The result indicates that the optimal total cost of the GSHP/ASHP combination system is $19.1 million for a 20-year operation, where 88% of the cooling load and 98% of the heating load can be covered by GSHP with 231*200m boreholes. These 231 boreholes can fit well into the empty area of the node by assuming a 6 m spacing. A Modelica model of the hybrid GSHP system was built for the detailed dynamic case study of the node. The correctness of the model was verified. It was found from the case study results that an appropriate hybrid configuration can not only help decrease the installation cost (the required borehole number) by decreasing the load to the ground but also help decrease the operation cost (the electricity consumption) by improving the efficiency of GSHP.